# ON PEIXOTO'S CONJECTURE FOR FLOWS ON NON-ORIENTABLE 2-MANIFOLDS

CARLOS GUTIERREZ AND BENITO PIRES

ABSTRACT. Contrary to the case of vector fields on orientable compact 2-manifolds, there is a smooth vector field $X$ on a non-orientable compact 2-manifold with a dense orbit (and therefore without closed orbits) whose phase portrait –up to topological equivalence– remains intact under a one-parameter family of twist perturbations localized in a flow box of $X$.

2000 *Mathematics Subject Classification.* Primary 34D30, 37E05, 37E35; Secondary 37C20.
The first author was supported in part by Pronex/CNPq/MCT grant number 66.2249/1997-6.
The second author was supported by Fapesp grant number 01/04598-0.

Departamento de Matemática, Institituo de Ciências Matemáticas e de Computaçâo, Universidade de São Paulo, Av. do Trabalhador São Carlense, 400, Centro, CEP 13560-970 São Carlos - SP, Brazil

*E-mail address*: gutp@icmc.usp.br

Departamento de Matemática, Institituo de Ciências Matemáticas e de Computaçâo, Universidade de São Paulo, Av. do Trabalhador São Carlense, 400, Centro, CEP 13560-970 São Carlos - SP, Brazil

*E-mail address*: bpires@icmc.usp.br